\documentclass[a4paper,11pt]{amsart}

\usepackage{amsthm}
\usepackage[english]{babel}
\usepackage{amsmath}
\usepackage{amssymb}
\usepackage{german}
\usepackage{hyperref}
\usepackage{booktabs}
\usepackage{enumerate}
\usepackage[left = 25mm, top = 40mm, bottom = 35mm, right = 25mm]{geometry}
\usepackage{multicol}
\usepackage{bbm}
\usepackage{todonotes}
\usepackage{setspace}
\usepackage{mathtools}
\usepackage {url}
\usepackage[capitalise]{cleveref}
\usepackage{xcolor}

\usepackage{tikz}

\makeatletter
\def\paragraph{\@startsection{paragraph}{4}%
  \z@\z@{-\fontdimen2\font}%
  {\normalfont\bfseries}}
\makeatother

\theoremstyle{definition}
\newtheorem{definition}{Definition}
\newtheorem{remark}{Remark}
\newtheorem{example}{Example}
\theoremstyle{plain}
\newtheorem{theorem}{Theorem}

\newtheorem{proposition}{Proposition}
\newtheorem{corollary}{Corollary}

\newtheorem*{question*}{Question}
\newtheorem*{theorem*}{Theorem}

\definecolor{defblue}{rgb}{0.1,0.1,0.7}
\newcommand{\defi}[1]{\textcolor{defblue}{\emph{#1}}}

\newcommand{\Z}{\mathbb{Z}}
\newcommand{\R}{\mathbb{R}}
\newcommand{\cM}{\mathcal{M}}
\newcommand{\cC}{\mathcal{C}}

\DeclareMathOperator{\Dcup}{\dot{\bigcup}}
\DeclareMathOperator{\bx}{\square}
\DeclareMathOperator{\smallsum}{\Sigma}

\newcommand{\set}[2]{\left\{ #1 \;\middle|\; #2 \right\}} 
\newcommand{\eqdef}{\mbox{\,\raisebox{0.2ex}{\scriptsize\ensuremath{\mathrm:}}\ensuremath{=}\,}} 
\newcommand{\defeq}{\mbox{~\ensuremath{=}\raisebox{0.2ex}{\scriptsize\ensuremath{\mathrm:}} }} 
\newcommand{\ssm}{\smallsetminus} 
\DeclareMathOperator{\inv}{inv} 
\DeclareMathOperator{\asc}{asc} 
\newcommand{\OEIS}[1]{{\rm \href{http://oeis.org/#1}{\texttt{#1}}}}

\title{Minimum maximal matchings in permutahedra}

\newcommand{\perm}{\Pi}

\author[S.~Brenner]{Sofia Brenner}
\address[S.~Brenner]{Institut f\"ur Mathematik, Universit\"at Kassel, Germany}
\email{sofia.brenner@mathematik.uni-kassel.de}
\urladdr{\url{https://sofiabrenner.github.io/}}

\author[J.~Fink]{Ji\v{r}í Fink}
\address[J.~Fink]{Charles University Prague, Czech Republic}
\email{jiri.fink@matfyz.cuni.cz}
\urladdr{\url{https://ktiml.mff.cuni.cz/~fink/}}

\author[H.~P.~Hoang]{Hung. P. Hoang}
\address[H.~P.~Hoang]{Algorithm and Complexity Group, Faculty of Informatics, TU Wien, Austria}
\email{phoang@ac.tuwien.ac.at}
\urladdr{\url{https://people.inf.ethz.ch/hoangp/}}

\author[A.~Merino]{Arturo Merino}
\address[A.~Merino]{Instituto de Ciencias de la Ingeniería, Universidad de O'Higgins, Chile}
\email{arturo.merino@uoh.cl}
\urladdr{\url{https://amerino.cl}}

\author[V.~Pilaud]{Vincent Pilaud}
\address[V.~Pilaud]{Universitat de Barcelona \& Centre de Recerca Matemàtica, Barcelona, Spain}
\email{vincent.pilaud@ub.edu}
\urladdr{\url{https://www.ub.edu/comb/vincentpilaud/}}

\date{\today}

\thanks{
This work was initiated at the 3rd Combinatorics, Algorithms, and Geometry workshop in Dresden, Germany in 2024.
We would like to thank the organizers and all the participants of the workshop for the inspiring atmosphere.
Especially, we thank L\'aszl\'o Kozma for stimulating discussions on this paper.
This work is also supported by Czech Science Foundation grant GA~22-15272S. \\
\indent We are indebted to Nathan Carter for providing the 3d coordinates of the embedding of \cref{fig:matching5} which can be originally found at \url{https://groupexplorer.sourceforge.net/images/cd-s5-transpo.gif}. \\
\indent Sofia Brenner received funding from the European Research Council (ERC) under the European Union’s Horizon 2020 research and innovation programme (EngageS: grant agreement No.~820148) and from the German Research Foundation DFG (SFB-TRR 195 ``Symbolic Tools in Mathematics and their Application'' as well as grant 522790373). \\
\indent Hung P. Hoang acknowledges support from the Austrian Science Foundation (FWF, project Y1329 START-Programm). \\
\indent Vincent Pilaud was partially supported by the Spanish project PID2022-137283NB-C21 of MCIN/AEI/10.13039/501100011033 / FEDER, UE, by the Severo Ochoa and María de Maeztu Program for Centers and Units of Excellence in R\&D (CEX2020-001084-M), by the Departament de Recerca i Universitats de la Generalitat de Catalunya (2021 SGR 00697), by the French project CHARMS (ANR-19-CE40-0017), and by the French--Austrian project PAGCAP (ANR-21-CE48-0020 \& FWF I 5788).
}

\begin{document}
	
\maketitle

\begin{abstract}
	We prove that the minimal size~$\cM(\perm_n)$ of a maximal matching in the permutahedron~$\perm_n$ is asymptotically~$n!/3$.
	On the one hand, we obtain a lower bound $\cM(\perm_n) \ge n! (n-1) / (3n-2)$ by considering $4$-cycles in the permutahedron.
	On the other hand, we obtain an asymptotical upper bound~$\cM(\perm_n) \le n!(1/3+o(1))$ by multiple applications of Hall's theorem (similar to the approach of Forcade for the hypercube~\cite{MR0321804}) and an exact upper bound~$\cM(\perm_n) \le n!/3$ by an explicit construction.
	We also derive bounds on minimum maximal matchings in products of permutahedra.
	
\end{abstract}

\section{Introduction}
Matchings are a fundamental concept in mathematics and computer science. 
Different variants of matchings can be used to model problems or as subroutines in algorithms.
While perfect matchings and maximum cardinality matchings are well understood, the landscape for small maximal matchings is less clear.
More specifically, given a graph $G$, we are interested in 
\[
 \cM (G) \eqdef \min \set{ |M| }{ M \text{ is a maximal matching in } G}.
\]
(The quantity $\cM(G)$ is also known as the \defi{edge domination number} of~$G$~\cite{MR579424}.)

For example, while a maximum cardinality matching can be found in polynomial time~\cite{MR177907}, finding a minimum maximal matching is an NP-hard problem~\cite{MR579424}, even when the graph is regular and bipartite~\cite{MR2472701}.
Assuming the Unique Games Conjecture, it is also NP-hard to approximate a minimum maximal matching with a constant better than two~\cite{MR3950886}. 
For general bounds, there are a couple of results.
First, $\cM(G) \geq \frac{m}{2\Delta - 1}$ where $m$ is the number of edges and $\Delta$ is the maximum degree~\cite{MR2074836}.
Second, $\frac{m}{\Delta^L} \leq \cM(G) \leq m - \Delta^L$ where $m$ is the number of edges and  $\Delta^L$ is the maximum degree of the line graph of $G$~\cite{MR4569919}. 
Thus, it is natural to ask for estimates and solutions for this problem for special graph~classes.

Let~$Q_n$ be the graph with $2^n$ vertices corresponding to the binary strings of length $n$, and where two vertices are adjacent if their corresponding strings differ in exactly one position.
It is the skeleton of the $n$-dimensional \defi{hypercube}, and the cover graph of the \defi{Boolean lattice}.
Let~$\perm_n$ be the graph with $n!$ vertices corresponding to the permutations of $[n]$, and where two vertices are adjacent if they differ by a transposition of two adjacent elements. 
It is the skeleton of the $(n-1)$-dimensional \defi{permutahedron}, and the cover graph of the \defi{weak order}.

Concerning the minimal size of a maximal matching in the hypercube~$Q_n$, Forcade~\cite{MR0321804} showed~that 
\[
\lim_{n\to \infty} \frac{1}{2^n}\cM(Q_n) = \frac{1}{3}.
\]
Our main contribution is an analogous result for the permutahedron $\perm_n$.
We obtain the following asymptotically tight bounds for $\cM(\perm_n)$.

\begin{theorem}
\label{thm:perm_tight}
The minimal size~$\cM(\perm_n)$ of a maximal matching of the permutahedron~$\perm_n$ is bounded by
\[
\frac{n-1}{3n-2}n!  \leq \cM(\perm_n) \leq \frac{1}{3} n!.
\]
\end{theorem}

The lower bound is obtained from a more general bound, where we argue in terms of the distribution of 4-cycles in any graph (\cref{sec:lower}).
This generalizes Forcade's lower bound argument for the hypercube~\cite{MR0321804} and is better than the aforementioned lower bounds.

The upper bound is obtained by a simple and explicit construction, obtained by combining maximal matchings of~$\perm_4$ in subgraphs of~$\perm_n$ (\cref{sec:upperboundpermutahedra}).
For any vertex, it takes linear time in $n$ to output its neighbor in the matching, if any.

We also discuss how Forcade's argument for the upper bound for the hypercube~\cite{MR0321804} can be packaged into a general framework for certain bipartite graphs (\cref{sec:upper_bipartite}), based on Hall's theorem.
Note that our general lower bound above works best on some regular graphs, while this framework for the upper bound only works on some bipartite graphs.
Hence, as a next step, we focus on the bipartite regular graphs for which we can obtain asymptotically tight bounds.

For this, we turn to the Cartesian products of permutahedra.
These objects also correspond to a rich class of polytopes, namely the bipartite regular quotientopes~\cite{hung_thesis}.
While the general bounds above are asymptotically tight in this case, we additionally give an explicit and tighter upper bound construction, based on our construction for the permutahedron (\cref{sec:cartesian}).

\begin{theorem}
\label{thm:quot_tight}
	Let $n_1, \dots, n_k \geq 2$ be integers, and let $n \eqdef n_1 + \dots + n_k$.
	Let $\perm$ be the Cartesian product of $\perm_{n_1}, \dots, \perm_{n_k}$. Then 
	\[
		\frac{n-k}{3n-3k+1}|V(\perm)| \leq \cM(\perm) \leq
		\begin{cases}
			\big(\frac{1}{3} + O(n^{-1/2})\big) |V(\perm)| & \text{if $\perm$ is a hypercube,} \\
			\frac{1}{3}|V(\perm)| & \text{otherwise.}
		\end{cases}
	\]
\end{theorem}

\section{Notation and preliminaries}

In this section, we introduce the notation used throughout this paper.
We denote the disjoint union by~$\sqcup$.

\paragraph{Matchings.}
A \defi{matching} $M$ in a graph $G \eqdef (V,E)$ is a subset of~$E$ such that every vertex in~$V$ is incident to at most one edge in $M$.
A vertex in~$V$ is \defi{covered} (resp.~\defi{exposed}) if it is incident to one (resp.~no) edge in~$M$. 
The matching~$M$ \defi{saturates} a subset~$X$ of~$V$ if all vertices in $X$ are covered.
A \defi{maximal matching} is a matching that is maximal with respect to inclusion.

\paragraph{Permutations.} 
Let $S_n$ denote the \defi{symmetric group}, that is, the group of permutations of~$[n] \eqdef$ $\{1, \dots, n\}$.
We write a permutation~$\sigma$ in \defi{one-line notation}, meaning as a word~$\sigma_1\dots \sigma_n$.
For a permutation $\sigma$, let $\inv(\sigma)$ denote its \defi{inversion set}, i.e., the set of pairs $(i, j)$ such that $i < j$ and $\sigma_i > \sigma_j$.
For $i \in [n-1]$, let $\tau_i \eqdef (i \; i+1)$ be the \defi{simple transposition} that only exchanges~$i$ and~$i+1$. 
We obtain $\sigma\tau_i$ from $\sigma$ by swapping the entries at positions $i$ and $i+1$. 

\paragraph{Permutahedra.} 
The \defi{permutahedron} $\perm_n$ is the Cayley graph of~$S_n$ for the simple transpositions~$\{\tau_1, \dots, \tau_{n-1}\}$.
In other words, its vertices are the permutations of $[n]$ and its edges are the pairs of permutations $\{\sigma, \sigma'\}$ such that $\sigma = \sigma' \tau_i$ for some $i \in [n-1]$. 
In particular, $\perm_n$ is $(n-1)$-regular. Moreover, every edge in $\perm_n$ corresponds to a unique transposition $\tau_i$ for $i \in [n-1]$, and we call such an edge a \defi{$\tau_i$-edge}.
For $i \ne j$ in~$[n-1]$, we call~\defi{$\tau_i\tau_j$-cycle} any cycle obtained by alternating between $\tau_i$-edges and~$\tau_j$-edges.
The~$\tau_i\tau_j$-cycles are $4$-cycles if~$|i-j| > 1$ and $6$-cycles otherwise.


\section{Lower bound for general graphs}
\label{sec:lower}

In this section, we derive a lower bound for the size of a maximal matching in a graph $G$ in terms of the number of its 4-cycles (\cref{prop:lower_bound}).
Our bound specializes to that of~\cite{MR0321804} for the hypercube (\cref{ex:lb_hypercube}).
Applied to the permutahedron, this yields an asymptotically tight lower bound (\cref{coro:lb_permutahedron}).

\begin{definition}
	For $\alpha \in \Z_{> 0}$, we say that a graph~$G$ is \defi{$\alpha$-heavy} if, 
	for every edge~$e$ in $G$, there are at least $\alpha$ induced 4-cycles such that $e$ is the only common edge of any two of these cycles.
\end{definition}

\begin{proposition}
	\label{prop:lower_bound}
	If $G=(V,E)$ is $\alpha$-heavy and has average degree $d$ and maximum degree $\Delta$, then any maximal matching of~$G$ has cardinality at least
	\[
	\frac{d}{4\Delta-\alpha-2} |V|.
	\]
\end{proposition}

\begin{proof}
	Let $M$ be a maximal matching in $G$.
	An \defi{$M$-edge} is an edge in~$M$.
	\mbox{A \defi{$1$-edge}} (resp.~\defi{$2$-edge}) is an edge in $G$ that is incident to exactly one $M$-edge (resp.~two $M$-edges). 
	Note in particular that $1$- and $2$-edges are not in~$M$ since they are adjacent to some $M$-edges.
	We denote by $m_0$, $m_1$ and $m_2$ the number of \mbox{$M$-edges}, $1$-edges, and $2$-edges, respectively.

	Since $G$ has average degree~$d$ and all edges are either $M$-edges, $1$-edges, or $2$-edges, we obtain
	\begin{equation}
		d|V| = 2(m_0 + m_1 + m_2).
		\label{eq:a}
	\end{equation}

	Double counting the number of adjacent pairs of edges with precisely one $M$-edge, we obtain
	\begin{equation}
		m_1 + 2m_2 \le 2(\Delta-1)m_0,
		\label{eq:b}	
	\end{equation}
	since each $1$-edge (resp.~$2$-edge) is adjacent to one (resp.~two) $M$-edges, while each $M$-edge is adjacent to at most $2(\Delta-1)$ edges that are either $1$- or $2$-edges.

	Finally, double counting the number of adjacent pairs of edges with one $M$-edge and one~$2$-edge, we obtain
	\begin{equation}
		\alpha m_0 \le 2 m_2.
		\label{eq:c}
	\end{equation}
	Indeed, each $2$-edge is adjacent to precisely two $M$-edges.
	Conversely, consider an $M$-edge~$e$ and a 4-cycle~$C$ that contains $e$.
	One of the two edges of~$C$ adjacent to~$e$ must be a $2$-edge, since otherwise we could add the remaining edge of~$C$ to~$M$ to create a larger matching, contradicting the maximality of $M$.
	The inequality~\eqref{eq:c} thus follows from the fact that~$G$ is $\alpha$-heavy.
	
	The sum \eqref{eq:a} + 2\eqref{eq:b} + \eqref{eq:c} gives $d|V| \le (4\Delta-\alpha-2)m_0$.
\end{proof}

\begin{example}
\label{ex:lb_hypercube}
	For the  $n$-dimensional hypercube~$Q_n$, we have~$d = \Delta = n$ and~${\alpha = n-1}$ so that~$m \geq |V| n / (3n-1)$, recovering the bound in~\cite{MR0321804}.
\end{example}

In general, the lower bound of \cref{prop:lower_bound} works best when the graph is regular (i.e., $d = \Delta$) and $\alpha = \Delta - c$, for some constant $c$, ideally $c = 1$ (as in the case of the hypercube).
Applying \cref{prop:lower_bound} to the permutahedron yields the following bound.

\begin{corollary}
\label{coro:lb_permutahedron}
	Every maximal matching of $\perm_n$ has at least $n! (n-1) / (3n-2)$ edges.
\end{corollary}

\begin{proof}
	This is a direct application of \cref{prop:lower_bound} since $\perm_n$ is $(n-1)$-regular and $(n-4)$-heavy.
	Indeed, any edge~$e$ of~$\perm_n$ is a $\tau_i$-edge for some~$i \in [n-1]$.
	For any~${j \in [n-1]}$ with~$|i-j| > 1$, let~$C_j$ be the $4$-cycle obtained by alternating $\tau_i$- and $\tau_j$-edges, starting with~$e$.
	Then $e$ is the only common edge of any two~$C_j$ and~$C_k$ with~$j \ne k$, and there are at least~$n-4$ such cycles given by the different~${j \in [n-1] \ssm \{i-1,i,i+1\}}$.
\end{proof}

\begin{remark}
\enlargethispage{.2cm}
Besides the hypercubes and the permutahedra, \cref{prop:lower_bound} can also be applied to the graphs of other classical polytopes:
\begin{description}
\item[Associahedra.] The \defi{associahedron} is the graph whose vertices are the binary trees with $n$ internal nodes and whose edges are tree rotations. It is regular with degree~$n-1$ and $(n-5)$-heavy, so that any maximal matching in the associahedron has at least~$\frac{(n-1)}{(3n-1)(n+1)}\binom{2n}{n}$ edges.
Using integer linear programming we calculated the minimum size of maximal matchings in the associahedron for~$n \le 6$; see Table \ref{tab:associahedron}.
Since exposed vertices form an independent set, we also included the maximal size of an independent set.
\item[Coxeter permutahedra.] A \defi{Coxeter permutahedron} is the convex hull of a generic point under the action of a finite Coxeter group~$W$~\cite{Humphreys,BjornerBrenti}. It is regular of degree~$n$ and~$(n-\delta)$-heavy, where~$n$ is the rank of~$W$ and~$\delta$ is the maximal degree of the Dynkin diagram of~$W$, hence any maximal matching in the Coxeter permutahedron has at least~$\frac{n}{3n+\delta-2}|W|$ edges. 
\end{description}
In contrast, it does not apply as such to all graphical zonotopes, since they can fail to be $\alpha$-heavy (some edges might appear in no $4$-cycle). 
For graphical zonotopes, we would need an improved version of \cref{prop:lower_bound} averaging the number of disjoint $4$-cycles containing an edge.
\end{remark}

\begin{table}[t]
\begin{center}
\begin{tabular}{@{\;}r@{\qquad}r@{\qquad}r@{\qquad}r@{\qquad}r@{\;}}\hline\rule{0 pt}{2.5ex}
Inner nodes $n$ & Vertices & Edges & Matching & Independent \\[.4ex]\hline\rule{0 pt}{3ex}
2 & 2 & 1 & 1 & 1 \\
3 & 5 & 5 & 2 & 2 \\
4 & 14 & 21 & 5 & 6 \\
5 & 42 & 84 & 14 & 16 \\
6 & 132 & 330 & 44 & 50 \\\hline
\end{tabular}
\end{center}
\caption{The sizes of minimum maximal matchings and maximum independent sets in the associahedron for $n \le 6$.}
\label{tab:associahedron}
\end{table}


\section{Upper bound for bipartite graphs}
\label{sec:upper_bipartite}

\enlargethispage{.4cm}
In this section, we obtain an upper bound on the minimal size of a maximal matching of certain bipartite graphs (\cref{prop:ub_Hall}).
Our approach is similar to that of~\cite{MR0321804} for the hypercube (\cref{ex:ub_hypercube}).
Applied to the permutahedron, this yields an asymptotically tight upper bound (\cref{coro:ub_permutahedron}).
The proof of \cref{prop:ub_Hall} uses Hall's classical matching theorem.

\begin{theorem}[Hall's theorem~\cite{MR1581694}]
\label{thm:hall}
	Let $G$ be a bipartite graph with two parts $X$ and $Y$.
	Then $G$ has a matching that saturates $X$ if and only if for every $X' \subseteq X$, the number of neighbors of~$X'$ in $Y$ is at least $|X'|$.
\end{theorem}

\begin{proposition}
\label{prop:ub_Hall}
Consider a graph~$G \eqdef (V,E)$ such that
\begin{itemize}
\item[(i)] $V = V_0 \Dcup \dots \Dcup V_\ell$ and $E \cap (V_i \times V_j) \ne \varnothing$ implies~$|i-j| = 1$,
\item[(ii)] for any~$0 \le i < \lfloor \ell/2 \rfloor$ and~$X \subseteq V_i$, the number of neighbors of~$X$ in~$V_{i+1}$ is at least~$|X|$,
\item[(iii)] for any~$\lceil \ell/2 \rceil < i \le \ell$ and~$X \subseteq V_i$, the number of neighbors of~$X$ in~$V_{i-1}$ is as least~$|X|$.
\end{itemize}
Then~$G$ admits a maximal matching of cardinality at most
\[
|V|/3 + 6\max \{ |V_{\lfloor \ell/2 \rfloor}|, |V_{\lceil \ell/2 \rceil}| \}.
\]
\end{proposition}

\begin{proof}
For~$k \in [\ell-1]$, we write~$G_k^-$ (resp.~$G_k^+$, resp.~$G_k^\pm$) for the subgraph of~$G$ induced by~$V_{k-1} \cup V_k$ (resp.~$V_k \cup V_{k+1}$, resp.~$V_{k-1} \cup V_k \cup V_{k+1}$).
For~${0 < k < \lfloor \ell/2 \rfloor}$, we obtain by two applications of \cref{thm:hall} (Hall's theorem) using (ii) that
\begin{itemize}
\item there exists a matching~$M_k^-$ of~$G_k^-$ which saturates~$V_{k-1}$,
\item there exists a matching~$M_k^+$ of~$G_k^+$ so that a vertex of~$V_k$ is covered in~$M_k^+$ if and only if it is exposed in~$M_k^-$.
\end{itemize}
Therefore, the union~$M_k^\pm \eqdef M_k^- \cup M_k^+$ is a matching of~$G_k^\pm$ which covers both~$V_{k-1}$ and~$V_k$.
Moreover, $|M_k^\pm| = |V_k|$ since all edges of~$M_k^\pm$ are incident to~$V_k$.
Similarly, for $\lceil \ell/2 \rceil < k < \ell$, there is a matching~$M_k^\pm$ of~$G_k^\pm$ which covers both~$V_k$ and~$V_{k+1}$, and with~$|M_k^\pm| = |V_k|$.

Let~$p \in [3]$ be such that~$\sum_{k \in [\ell], \, k \equiv p \pmod{3}} |V_k| \le |V|/3$.
Let~$M$ be the union of the matchings~$M_k^\pm$ for~$k \in [\ell-1] \ssm \{ \lfloor \ell/2 \rfloor, \lceil \ell/2 \rceil \}$ with~$k \equiv p \pmod{3}$.
Note that~$M$ is a matching (since the graphs~$G_k^\pm$ for~$k \equiv p \pmod{3}$ are vertex disjoint) and that~$|M| \le |V|/3$ (since~$|M_k^\pm| = |V_k|$).
Let~$M'$ be a maximal matching containing~$M$.
Since~$M_k^\pm$ covers $V_{k-1}$ and~$V_k$ when~$0 < k < \lfloor \ell/2 \rfloor$ (resp.~$V_k$ and~$V_{k+1}$ when~$\lceil \ell/2 \rceil < k < \ell$), all edges of~$M' \ssm M$ are incident to~$V_0 \cup V_{\lfloor \ell/2 \rfloor} \cup V_{\lfloor \ell/2 \rfloor -1} \cup V_{\lceil \ell/2 \rceil} \cup V_{\lceil \ell/2 \rceil+1} \cup V_\ell$.
As~(ii) and~(iii) imply that~${|V_k| \le \max \{ |V_{\lfloor \ell/2 \rfloor}|, |V_{\lceil \ell/2 \rceil}| \}}$ for all~$k$, we indeed obtained a maximal matching~$M'$ with~$|M'| \le |V|/3 + 6 \max \{ |V_{\lfloor \ell/2 \rfloor}|, |V_{\lceil \ell/2 \rceil}| \}$.
\end{proof}

\begin{example}
\label{ex:ub_hypercube}
The $n$-dimensional hypercube $Q_n$ satisfies the assumptions of \cref{prop:ub_Hall}.
For~$0 \le k \le n$, denote by~$V_k$ the set of binary strings of length $n$ with precisely $k$ occurrences of~$1$.
If~$X \subseteq V_k$ and~$Y$ denotes its neighborhood in~$V_{k+1}$, we have~$(n-k) |X| = (k+1) |Y|$, which proves~(ii) and~(iii) in \cref{prop:ub_Hall}.
It follows that $Q_n$ admits a maximal matching of cardinality at most~$2^n/3 + 6\binom{n}{\lfloor n/2 \rfloor} = 2^n/3 \big( 1 +O(n^{-1/2}) \big)$.
This was precisely the approach of~\cite{MR0321804}.
\end{example}

Finally, we observe that \cref{prop:ub_Hall} yields an asymptotically tight bound for the minimal size of a maximal matching in the permutahedron.

\begin{corollary}
\label{coro:ub_permutahedron}
The permutahedron $\perm_n$ satisfies the assumptions of \cref{prop:ub_Hall}.
Hence, it admits a maximal matching of cardinality at most
\[
n! \Big( \frac{1}{3} + \frac{6}{\sqrt{2\pi}}n^{-3/2}  + o(n^{-3/2}) \Big).
\]
\end{corollary}

\begin{proof}
Recall that we denote by~$\inv(\sigma) \eqdef \set{(i,j) \in [n]^2}{i < j \text{ and } \sigma_i > \sigma_j}$ the inversion set and by~$\asc(\sigma) \eqdef \set{i \in [n-1]}{\sigma_i < \sigma_{i+1}}$ the ascent set of a permutation~$\sigma$ of~$[n]$.
For~$0 \le k \le \binom{n}{2} \defeq \ell$, we denote by~$V_i \eqdef \set{\sigma \in S_n}{|\inv(\sigma)| = k}$ the set of permutations with precisely~$k$ inversions.
Let~$\R S_n$ be the real vector space with basis indexed by~$S_n$.
Consider the linear map $U: \R S_n \to \R S_n$ defined by $U(\sigma) = \sum_{i \in \asc(\sigma)} i \cdot \sigma\tau_i$.
Gaetz and Gao~\cite{MR4042823} showed that $U^{\ell-2k}$ is an isomorphism between $\R V_k$ and $\R V_{\ell-k}$ for~$0 \le k \le \lfloor k/2 \rfloor$.
This immediately implies that~$U$ is injective.

Consider now a subset~$X$ of~$V_k$ and its neighborhood~$Y$ in~$V_{k+1}$.
Since~$U$ is injective, the image by~$U$ of the subspace of~$\R V_i$ generated by~$X$ has dimension~$|X|$.
But by definition of~$U$, this image is contained in the subspace of~$\R V_{k+1}$ generated by~$Y$, which has dimension~$|Y|$. 
Hence, we obtain that~${|X| \leq |Y|}$.
This shows~(ii), and the proof of~(iii) is symmetric.

Finally, $M_n \eqdef \max \{ |V_{\lfloor \ell/2 \rfloor}|, |V_{\lceil \ell/2 \rceil}| \}$ is known as the \emph{Kendall-Mann number} (\OEIS{A000140}), the row maximum of table of Mahonian numbers (\OEIS{A008302}).
It is known (it follows for instance from~\cite{MR2324533}) that~${M_n \sim 6 n^{n-1} / e^n \sim n! \cdot n^{-3/2} \cdot 6/\sqrt{2\pi}}$, which yields our asymptotic bound.
\end{proof}

\begin{remark}
\label{rem:scd}
By using Hall's theorem, the proof of \cref{prop:ub_Hall} is not constructive.
For the application to the hypercube in \cref{ex:ub_hypercube}, we can make it constructive by using a symmetric chain decomposition in the hypercube.
(Note that this decomposition is a well known concept in partial order theory; we define it here in graph theoretic terms.)
A \defi{symmetric chain} in the $n$-dimensional hypercube is a path $(x_t, x_{t+1}, \dots, x_{n-t})$ for some $t \in \lfloor n/2 \rfloor$ such that $x_i$ has exactly $i$ occurrences of 1 for $i\in\{t,\dots,n-t\}$.
A \defi{symmetric chain decomposition} in the hypercube is a partition of its vertices into symmetric chains.
Greene and Kleitman~\cite{MR389608} described a simple construction of such a decomposition as follows.
For a binary word~$x$, we interpret the 0s in $x$ as opening brackets and the 1s as closing brackets.
For example, if $x = 1000110$, we interpret $x$ as $)((())($.
By matching the opening and closing brackets in the natural way, we can obtain the symmetric chain containing $x$ by flipping the rightmost unmatched 1 or the leftmost unmatched 0, until no more unmatched bits can be flipped.
In the example above, the chain that contains $x$ is $(0000110, 1000110, 1100110, 1100111)$.
The chains in this decomposition can be used to explicitly describe the matchings $M_k^+$ and $M_k^-$ in the proof of \cref{prop:ub_Hall} for the hypercube.
For the permutahedron, however, it is unknown whether an analogous symmetric chain decomposition exists, and as far as we are aware, there is no constructive proof of~(ii) and~(iii).
\end{remark}


\section{Explicit construction for permutahedra}
\label{sec:upperboundpermutahedra}

In this section, we construct an explicit maximal matching of size $n!/3$ in the permutahedron~$\perm_n$.
This bound is constructive and tighter than the upper bound in \cref{coro:ub_permutahedron}, and thus gives an alternative proof that the lower bound of~\cref{coro:lb_permutahedron} is asymptotically tight.
Our construction involves three steps: 
\begin{enumerate}
\item We first define two maximal matchings~$M_+$ and~$M_-$ of~$\perm_4$ with $8$ edges.
\item We then construct a matching~$M$ of~$\perm_n$ with~$n!/3$ edges by transporting~$M_+$ and~$M_-$ to certain maximal matchings in each subgraph of~$\perm_n$ induced by permutations with a fixed suffix of length $n-4$.
\item We finally prove that~$M$ is maximal.
\end{enumerate}

\subsection{Two maximal matchings of~$\perm_4$}

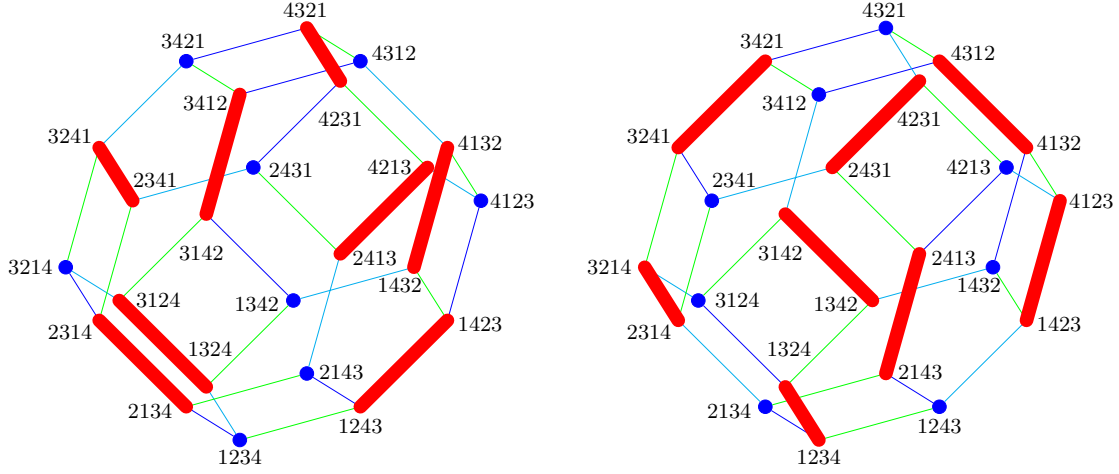
\begin{figure}[t]
\centerline{\scalebox{.9}{
\begin{tikzpicture}[scale = 1.8, font=\footnotesize]
	\coordinate (4132) at (1.414, 0.707,0) {};
	\coordinate (1432) at (1.414, 0,0.707) {};
	\coordinate (1423) at (1.414, -0.707,0) {};
	\node[draw, circle, fill=blue, inner sep=2pt, blue] (4123) at (1.414, 0,-0.707) {};	
	\coordinate (3241) at (-1.414, 0.707,0) {};
	\node[draw, circle, fill=blue, inner sep=2pt, blue] (3214) at (-1.414, 0,0.707) {};
	\coordinate (2314) at (-1.414, -0.707,0) {};
	\coordinate (2341) at (-1.414, 0,-0.707) {};
	\node[draw, circle, fill=blue, inner sep=2pt, blue] (4312) at (0.707, 1.414,0) {};
	\coordinate (3412) at (0,1.414,0.707) {};
	\node[draw, circle, fill=blue, inner sep=2pt, blue] (3421) at (-0.707, 1.414, 0) {};
	\coordinate (4321) at (0,1.414, -0.707) {};
	\coordinate (1243) at (0.707, -1.414,0) {};
	\node[draw, circle, fill=blue, inner sep=2pt, blue] (1234) at (0,-1.414,0.707) {};
	\coordinate (2134) at (-0.707, -1.414, 0) {};
	\node[draw, circle, fill=blue, inner sep=2pt, blue] (2143) at (0,-1.414, -0.707) {};
	\node[draw, circle, fill=blue, inner sep=2pt, blue] (1342) at (0.707, 0, 1.414) {};
	\coordinate (3142) at (0,0.707, 1.414) {};
	\coordinate (3124) at (-0.707, 0,1.414) {};
	\coordinate (1324) at (0,-0.707,1.414) {};
	\coordinate (4213) at (0.707, 0, -1.414) {};
	\coordinate (4231) at (0,0.707, -1.414) {};
	\node[draw, circle, fill=blue, inner sep=2pt, blue] (2431) at (-0.707, 0, -1.414) {};
	\coordinate (2413) at (0,-0.707,-1.414) {};
	
	\node[] (L4132) at (1.68, 0.74,0) {4132};
	\node[] (L1432) at (1.3, -0.15,0.707) {1432};	
	\node[] (L1423) at (1.68, -0.74,0) {1423};		
	\node[] (L4123) at (1.67, 0,-0.707) {4123};
	\node[] (L3241) at (-1.65, 0.8,0) {3241};
	\node[] (L3214) at (-1.7, 0,0.707) {3214};
	\node[] (L2314) at (-1.65, -0.8,0) {2314};	
	\node[] (L2341) at (-1.23, 0.17, -0.707) {2341};
	\node[] (L4312) at (0.98, 1.48, 0) {4312};
	\node[] (L1234) at (0,-1.55,0.707) {1234};
	\node[] (L3421) at (-0.73, 1.58, 0) {3421};
	\node[] (L2143) at (0.28,-1.4, -0.707) {2143};
	\node[] (L1243) at (0.707, -1.57,0) {1243};
	\node[] (L3412) at (-0.28,1.34,0.707) {3412};
	\node[] (L2134) at (-1, -1.45, 0) {2134};	
	\node[] (L4321) at (-0.05, 1.52, -0.807) {4321};
	\node[] (L1342) at (0.41, -0.03, 1.414) {1342};
	\node[] (L3142) at (-0.024, 0.4, 1.454) {3142};	
	\node[] (L3124) at (-0.39, 0,1.414) {3124};	
	\node[] (L1324) at (0.03,-0.4,1.414) {1324};
	\node[] (L4213) at (0.4, 0, -1.414) {4213};
	\node[] (L4231) at (0,0.4, -1.414) {4231};
	\node[] (L2431) at (-0.4, -0.02, -1.414) {2431};
	\node[] (L2413) at (0.28,-0.74,-1.414) {2413};
	
	\draw[blue] 
	(4132) -- (1432) 
	(1423) -- (4123)
	(3214) -- (2314)
	(2341) -- (3241)
	(4312) -- (3412)
	(3421) -- (4321)
	(1234) -- (2134)
	(2143) -- (1243)
	(1342) -- (3142)
	(3124) -- (1324)
	(4231) -- (2431)
	(2413) -- (4213)
	;

	\draw[cyan] 
	(4132) -- (4312)
	(4123) -- (4213)
	(2143) -- (2413)
	(2341) -- (2431)
	(4321) -- (4231)
	(1432) -- (1342)
	(3412) -- (3142)
	(1423) -- (1243)
	(3124) -- (3214)
	(1324) -- (1234)
	(3421) -- (3241)
	(2314) -- (2134)
	;

	\draw[green] 
	(1432) -- (1423)
	(4123) -- (4132)
	(3241) -- (3214)
	(2314) -- (2341)
	(3412) -- (3421)
	(4321) -- (4312)
	(1243) -- (1234)
	(2134) -- (2143)
	(3142) -- (3124)
	(1324) -- (1342)
	(4213) -- (4231)
	(2431) -- (2413)
	;
	
	\draw[line width= 2mm, red, line cap=round] 
	(2134) -- (2314)
	(3241) -- (2341)
	(3412) -- (3142)
	(2413) -- (4213)
	(4132) -- (1432)
	(1243) -- (1423)
	(3124) -- (1324)
	(4321) -- (4231)
	;
\end{tikzpicture}
\quad
\begin{tikzpicture}[scale = 1.8, font=\footnotesize]
	\coordinate (4132) at (1.414, 0.707,0) {};
	\node[draw, circle, fill=blue, inner sep=2pt, blue] (1432) at (1.414, 0,0.707) {};
	\coordinate (1423) at (1.414, -0.707,0) {};
	\coordinate (4123) at (1.414, 0,-0.707) {};
	\coordinate (3241) at (-1.414, 0.707,0) {};
	\coordinate (3214) at (-1.414, 0,0.707) {};
	\coordinate (2314) at (-1.414, -0.707,0) {};
	\node[draw, circle, fill=blue, inner sep=2pt, blue] (2341) at (-1.414, 0,-0.707) {};
	\coordinate (4312) at (0.707, 1.414,0) {};
	\node[draw, circle, fill=blue, inner sep=2pt, blue] (3412) at (0,1.414,0.707) {};
	\coordinate (3421) at (-0.707, 1.414, 0) {};
	\node[draw, circle, fill=blue, inner sep=2pt, blue] (4321) at (0,1.414, -0.707) {};
	\node[draw, circle, fill=blue, inner sep=2pt, blue] (1243) at (0.707, -1.414,0) {};
	\coordinate (1234) at (0,-1.414,0.707) {};
	\node[draw, circle, fill=blue, inner sep=2pt, blue] (2134) at (-0.707, -1.414, 0) {};
	\coordinate (2143) at (0,-1.414, -0.707) {};
	\coordinate (1342) at (0.707, 0, 1.414) {};
	\coordinate (3142) at (0,0.707, 1.414) {};
	\node[draw, circle, fill=blue, inner sep=2pt, blue] (3124) at (-0.707, 0,1.414) {};
	\coordinate (1324) at (0,-0.707,1.414) {};
	\node[draw, circle, fill=blue, inner sep=2pt, blue] (4213) at (0.707, 0, -1.414) {};
	\coordinate (4231) at (0,0.707, -1.414) {};
	\coordinate (2431) at (-0.707, 0, -1.414) {};
	\coordinate (2413) at (0,-0.707,-1.414) {};
	
	\node[] (L4132) at (1.68, 0.74,0) {4132};
	\node[] (L1432) at (1.3, -0.15,0.707) {1432};	
	\node[] (L1423) at (1.68, -0.74,0) {1423};		
	\node[] (L4123) at (1.67, 0,-0.707) {4123};
	\node[] (L3241) at (-1.65, 0.8,0) {3241};
	\node[] (L3214) at (-1.7, 0,0.707) {3214};
	\node[] (L2314) at (-1.65, -0.8,0) {2314};	
	\node[] (L2341) at (-1.23, 0.17, -0.707) {2341};
	\node[] (L4312) at (0.98, 1.48, 0) {4312};
	\node[] (L1234) at (0,-1.55,0.707) {1234};
	\node[] (L3421) at (-0.73, 1.58, 0) {3421};
	\node[] (L2143) at (0.28,-1.4, -0.707) {2143};
	\node[] (L1243) at (0.707, -1.57,0) {1243};
	\node[] (L3412) at (-0.28,1.34,0.707) {3412};
	\node[] (L2134) at (-1, -1.45, 0) {2134};	
	\node[] (L4321) at (-0.05, 1.52, -0.807) {4321};
	\node[] (L1342) at (0.41, -0.03, 1.414) {1342};
	\node[] (L3142) at (-0.024, 0.4, 1.454) {3142};	
	\node[] (L3124) at (-0.39, 0,1.414) {3124};	
	\node[] (L1324) at (0.03,-0.4,1.414) {1324};
	\node[] (L4213) at (0.4, 0, -1.414) {4213};
	\node[] (L4231) at (0,0.4, -1.414) {4231};
	\node[] (L2431) at (-0.4, -0.02, -1.414) {2431};
	\node[] (L2413) at (0.28,-0.74,-1.414) {2413};
	
	\draw[blue] 
	(4132) -- (1432) 
	(1423) -- (4123)
	(3214) -- (2314)
	(2341) -- (3241)
	(4312) -- (3412)
	(3421) -- (4321)
	(1234) -- (2134)
	(2143) -- (1243)
	(1342) -- (3142)
	(3124) -- (1324)
	(4231) -- (2431)
	(2413) -- (4213)
	;

	\draw[cyan] 
	(4132) -- (4312)
	(4123) -- (4213)
	(2143) -- (2413)
	(2341) -- (2431)
	(4321) -- (4231)
	(1432) -- (1342)
	(3412) -- (3142)
	(1423) -- (1243)
	(3124) -- (3214)
	(1324) -- (1234)
	(3421) -- (3241)
	(2314) -- (2134)
	;

	\draw[green] 
	(1432) -- (1423)
	(4123) -- (4132)
	(3241) -- (3214)
	(2314) -- (2341)
	(3412) -- (3421)
	(4321) -- (4312)
	(1243) -- (1234)
	(2134) -- (2143)
	(3142) -- (3124)
	(1324) -- (1342)
	(4213) -- (4231)
	(2431) -- (2413)
	;
	
	\draw[line width= 2mm, red, line cap=round] 
	(1324) -- (1234)
	(2143) -- (2413)
	(4231) -- (2431)
	(3214) -- (2314)
	(4123) -- (1423)
	(3142) -- (1342)
	(4312) -- (4132)
	(3241) -- (3421)
	;
\end{tikzpicture}
}}
\caption{Matching edges (red) and exposed vertices (blue) for $M^+$ (left) and $M^-$ (right). All other edges are colored according to the exchange position.}
\label{fig:matchings}
\end{figure}

We consider the two matchings $M^+$ and $M^-$ in $\perm_4$ of \cref{fig:matchings}.
Note~that
\begin{itemize}
\item both~$M^+$ and~$M^-$ are maximal matchings of~$\perm_4$ with $8$ edges,
\item the sets~$E^+$ and~$E^-$ of exposed vertices of~$M^+$ and $M^-$ are disjoint.
\end{itemize}

\subsection{Combining maximal matchings}

We now combine copies of the matchings~$M^+$ and~$M^-$ to create a matching~$M$ of~$\perm_n$ with~$n!/3$ edges.
Let~$S$ be the set of duplicate-free strings in $[n]$ of length~${n-4}$.
For~$s \in S$, denote by~$\perm^s$ the subgraph of~$\perm_n$ induced by permutations with suffix~$s$.
Note that~$\perm_n = \bigsqcup_{s \in S} \perm^s$.
Denote by~${\bar s = \{\bar s_1 < \bar s_2 < \bar s_3 < \bar s_4\}}$ the set of elements in $[n]$ that do not occur in~$s$.
For a permutation~$\pi$ of~$[4]$, let~$\bar s_\pi \eqdef \bar s_{\pi_1} \bar s_{\pi_2} \bar s_{\pi_3} \bar s_{\pi_4}$ and define~$\phi_s(\pi) \eqdef \bar  s_\pi s \in \perm^s$.
Observe~that $\phi_s$ defines a graph isomorphism from~$\perm_4$ to~$\perm^s$.
Define~$\varepsilon(s) \eqdef (-1)^{|\inv(s)| + \smallsum(\bar{s})}$, where $\smallsum(\bar s) \eqdef \bar s_1 + \bar s_2 + \bar s_3 + \bar s_4$ and~$\inv(s) \eqdef \set{(i,j) \in [n-4]^2}{i < j \text{ and } s_i > s_j}$ is the inversion set of~$s$.
Finally, define
\[
M \eqdef \bigsqcup_{s \in S} 
\phi_s(M^{\varepsilon(s)}).
\]

\begin{figure}[t]
\centerline{\includegraphics[width=\textwidth]{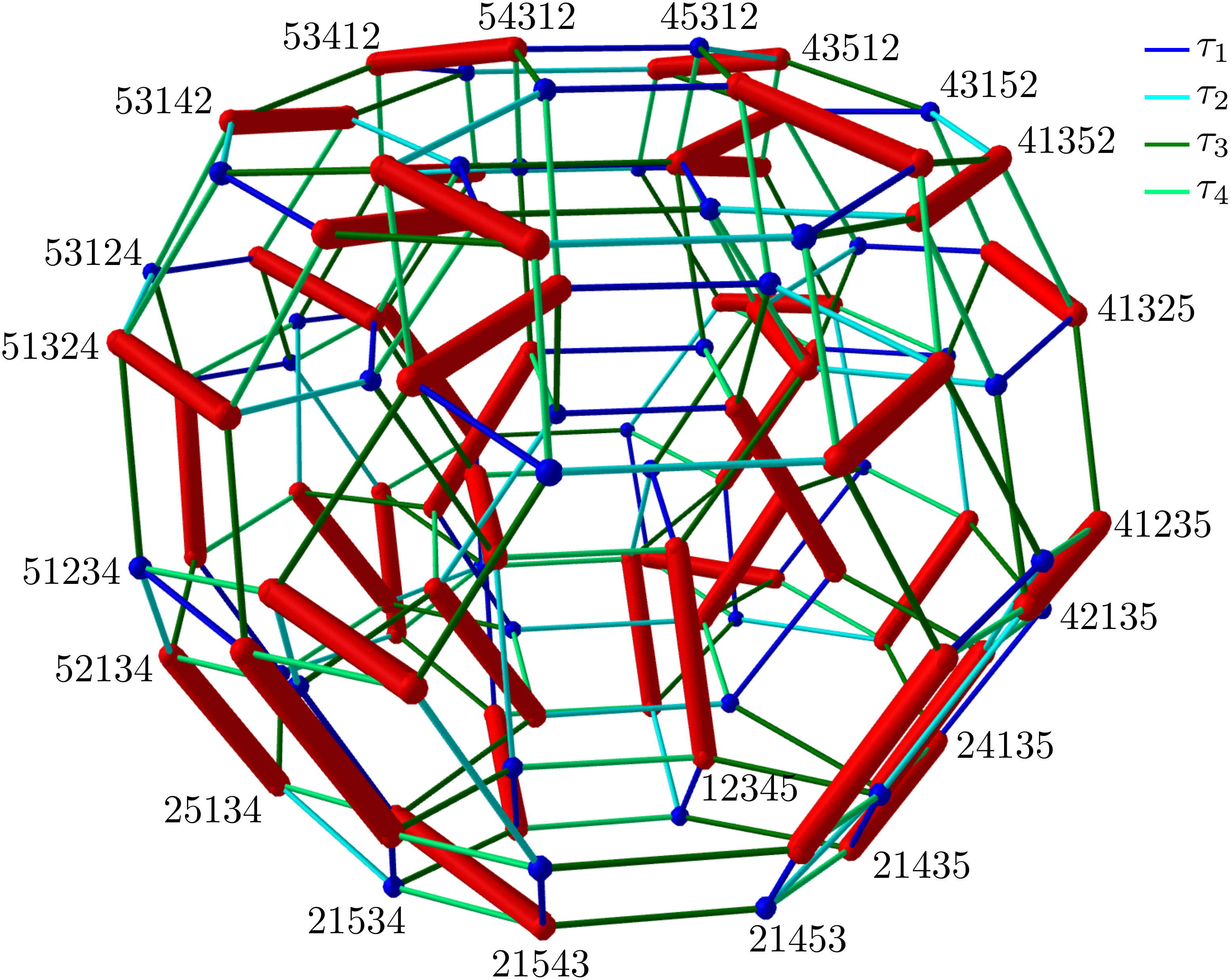}}
\caption{Matching edges (red) and exposed vertices (blue) for the matching~$M$ when~${n = 5}$. All other edges are colored according to the exchange position. The coordinates of the embedding were communicated to us by Nathan Carter. An animated 3d version can be found at \url{https://tinyurl.com/maximalMatchingPermutahedron}.}
\label{fig:matching5}
\end{figure}

\begin{example}
For instance, \cref{fig:matching5} illustrates the matching~$M$ when~$n = 5$.
Its edges are
\[
\begin{array}{c@{\quad}c@{\quad}c@{\quad}c@{\quad}c}
(23451, 24351) & (13542, 15342) & (12453, 14253) & (12534, 15234) & (12345, 13245) \\
(24531, 42531) & (15432, 51432) & (14523, 41523) & (15324, 51324) & (13425, 31425) \\ 
(25341, 52341) & (14352, 41352) & (15243, 51243) & (13254, 31254) & (14235, 41235) \\
(32541, 35241) & (31452, 34152) & (21543, 25143) & (21354, 23154) & (21435, 24135) \\
(34251, 43251) & (34512, 43512) & (24153, 42153) & (23514, 32514) & (23145, 32145) \\
(35421, 53421) & (35142, 53142) & (25413, 52413) & (25134, 52134) & (24315, 42315) \\
(43521, 45321) & (41532, 45132) & (42513, 45213) & (31524, 35124) & (32415, 34215) \\ 
(52431, 54231) & (53412, 54312) & (51423, 54123) & (52314, 53214) & (41325, 43125)
\end{array}
\]
and its exposed vertices are
\[
\begin{array}{c@{\;\;}c@{\;\;}c@{\;\;}c@{\;\;}c@{\;\;}c@{\;\;}c@{\;\;}c@{\;\;}c@{\;\;}c}
12354 & 13524 & 21345 & 23541 & 31245 & 34125 & 41253 & 43215 & 51234 & 53241 \\ 
12435 & 14325 & 21453 & 24513 & 31542 & 34521 & 42135 & 45123 & 51342 & 54132 \\ 
12543 & 14532 & 21534 & 25314 & 32154 & 35214 & 42351 & 45231 & 52143 & 54213 \\
13452 & 15423 & 23415 & 25431 & 32451 & 35412 & 43152 & 45312 & 53124 & 54321
\end{array}
\]
\end{example}

\begin{theorem}
\label{thm:ub_permutahedron}
The set $M$ is a maximal matching in $\perm_n$ of size $n!/3$. 
\end{theorem}

\begin{proof}
We have~$|M| = \sum_{s \in S} |M^{\varepsilon(s)}| = 8|S| = n!/3$ since~$M = \bigsqcup_{s \in S} \phi_s(M^{\varepsilon(s)})$.
Moreover, $M$ is a matching since~$\perm_n = \bigsqcup_{s \in S} \perm^s$ and any edge of~$\phi_s(M^{\varepsilon(s)})$ lies in~$\perm^s$.
We thus just need to prove that~$M$ is maximal, and this is the purpose of the next section.
\end{proof}

\begin{remark}
\label{rem:secondmatchingpn}
For later use, we point out that we can in fact define two maximal matchings~$M^\bullet \eqdef \bigsqcup_{s \in S} \phi_s(M^{\varepsilon(s)})$ and~$M^\circ \eqdef \bigsqcup_{s \in S} \phi_s(M^{-\varepsilon(s)})$ of the permutahedron~$\perm_n$ of size~$n!/3$.
The sets of exposed vertices~$E^\bullet$ and~$E^\circ$ of these matchings~$M^\bullet$ and~$M^\circ$ are given by
\[
E^\bullet = \bigsqcup_{s \in S} \phi_s(E^{\varepsilon(s)})
\qquad\text{and}\qquad
E^\circ = \bigsqcup_{s \in S} \phi_s(E^{-\varepsilon(s)}).
\]
Hence, $E^\bullet \cap E^\circ = \varnothing$ since~$E^+ \cap E^- = \varnothing$.
\end{remark}

\subsection{Proof of maximality}

In this section, we prove that~$M$ is maximal.
Assume by means of contradiction that~$M$ admits two exposed vertices~$\sigma, \sigma' \in S_n$ with~$\sigma = \sigma \tau_i$ for some~${i \in [n-1]}$.
Let~$s$ be the suffix formed by the last~$n-4$ letters of~$\sigma$, and~$\pi$ be the permutation of~$[4]$ such that~$\sigma = \phi_s(\pi)$.
Similarly, let~$s'$ and~$\pi'$ be such that~$\sigma' = \phi_{s'}(\pi')$.

If~$i < 4$, then~$s = s'$ so that~$\sigma$ and~$\sigma'$ belong to the same~$\perm^s = \perm^{s'}$. This contradicts the maximality of~$M^{\varepsilon(s)}$.
If~$i > 4$, then~$\bar s = \bar s'$ and~$|\inv(s)| = |\inv(s')| \pm 1$, so that~$\varepsilon(s) \ne \varepsilon(s')$. As~$\bar s \in E^{\varepsilon(s)}$ and~$\bar s' \in E^{\varepsilon(s')}$, this contradicts that~$E^+ \cap E^- = \varnothing$.

We can thus assume from now on that~$i = 4$.
Assume moreover without loss of generality that~$\sigma_4 < \sigma_5$ and set~$t \eqdef |\set{j \in [3]}{\sigma_4 < \sigma_j < \sigma_5}|$.
For~$j \in [3]$, we have~$\bar s_{\pi_j} = \sigma_j = \sigma'_j = \bar s'_{\pi'_j}$.
Hence, we have
\begin{equation}
\label{eq:relpi}
\pi_j =
\begin{cases}
	\pi'_j+1 & \text{ if } \sigma_4 < \sigma_j < \sigma_5, \\
	\pi'_j & \text{ otherwise}. \\
\end{cases}
\end{equation}
Observing~$E^+$ and~$E^-$ in \cref{fig:matchings}, we therefore obtain that~$(\pi, \pi')$ belongs to
\begin{itemize}
\item $\{(1\bar3\bar4\bar2, 1\bar2\bar3\bar4), (\bar4\bar31\bar2, \bar3\bar21\bar4)\}$ if $\varepsilon(s) = \varepsilon(s') = +$,
\item $\{(\bar2\bar34\bar1, \bar1\bar24\bar3), (4\bar3\bar2\bar1, 4\bar2\bar1\bar3)\}$ if $\varepsilon(s) = \varepsilon(s') = -$,
\item $\{(21\bar4\bar3, 21\bar3\bar4), (1\bar34\bar2, 1\bar24\bar3), (21\bar4\bar3, 21\bar3\bar4), (\bar412\bar3, \bar312\bar4), (4\bar31\bar2, 4\bar21\bar3), (34\bar2\bar1, 34\bar1\bar2)\}$ if $\varepsilon(s) = +$ and~$\varepsilon(s') = -$,
\item $\{(12\bar4\bar3, 12\bar3\bar4), (\bar2\bar3\bar4\bar1, \bar1\bar2\bar3\bar4), (\bar234\bar1, \bar134\bar2), (\bar421\bar3, \bar321\bar4), (\bar4\bar3\bar2\bar1, \bar3\bar2\bar1\bar4), (43\bar2\bar1, 43\bar1\bar2)\}$ if~$\varepsilon(s) = -$ and~$\varepsilon(s') = +$.
\end{itemize}
In each such pair~$(\pi, \pi')$, we have overlined the positions~$j \in [4]$ where~$\pi_j \ne \pi'_j$.
Hence, as~$t = |\set{j \in [3]}{\pi_j \ne \pi'_j}$| by \cref{eq:relpi}, we obtain from this case analysis that
\begin{equation}
\label{eq:even}
\varepsilon(s) = \varepsilon(s')
\iff
t \text{ is even.}
\end{equation}
Observe now that since~$\sigma = \sigma'\tau_4$ and~$\sigma_4 < \sigma_5$, we have
\[
{|\inv(s)| = |\inv(s')| + \sigma_5 - \sigma_4 - t - 1}
\qquad\text{and}\qquad
\smallsum(\bar s) = \smallsum(\bar s') + \sigma_4 - \sigma_5.
\]
Hence, as $\varepsilon(s)$ just records the parity of~$|\inv(s)| + \smallsum(\bar s)$, we obtain that
\begin{equation}
\label{eq:odd}
\varepsilon(s) = \varepsilon(s')
\iff
t \text{ is odd.}
\end{equation}
This concludes the proof since \eqref{eq:even} and~\eqref{eq:odd} contradict each other.


\section{Cartesian products of permutahedra}
\label{sec:cartesian}

In this section, we prove \cref{thm:quot_tight} concerning minimum maximal matchings in Cartesian products of permutahedra.
Let us first recall the definition.

\begin{definition}
The \defi{Cartesian product} $G \bx H$ of two graphs $G$ and $H$ is the graph with vertex set $V(G)\times V(H)$ and edge set $\set{(u,v)(u',v)}{uu' \in E(G), v\in V(H)} \cup \set{(u,v)(u,v')}{u \in V(G), vv'\in E(H)}$.
\end{definition}

We start with the lower bound.

\begin{proposition}
\label{prop:lb_cartesian}
Let $n_1, \dots, n_k$ be integers with $n_1 \geq \dots \geq n_k \geq 2$ and~${n_1 \geq 3}$.
Let $n \eqdef n_1 + \dots + n_k$ and $\perm \eqdef \perm_{n_1} \bx \dots \bx \perm_{n_k}$.
Then any maximal matching of~$\perm$ has cardinality at least 
\[
\frac{n-k}{3n-3k+1} V(\perm).
\]
\end{proposition}

\begin{proof}
First, $\perm$ is $(n - k)$-regular since each $\perm_{n_i}$ is $(n_i - 1)$-regular.
Next, we prove by induction on $i \in [k]$ that $G_{i} \eqdef \perm_{n_1} \bx \dots \bx \perm_{n_i}$ is $(N_i-i-3)$-heavy, for $N_i \eqdef n_1 + \dots+ n_i$.
For the base case $i=1$, this follows from the proof of \cref{coro:lb_permutahedron}.
For the inductive step $i \geq 2$, note that $G_{i} = G_{i-1} \bx \perm_{n_i}$.
Hence, we can write each vertex of $G_i$ as $(x,y)$, where $x$ and $y$ are vertices in $G_{i-1}$ and~$\perm_{n_i}$, respectively.
Let $e$ be an edge of $G_i$ between two vertices $(x,y)$ and $(x',y')$.
By the definition of a Cartesian product, either $x = x'$ or $y = y'$.

Firstly consider the case that $x = x'$.
By the proof of \cref{coro:lb_permutahedron}, there are $n_i-4$ cycles of $\perm_{n_i}$ such that $x'y'$ is the only common edge of any two such cycles.
Adding $x$ as a prefix to all vertices to these cycles, we obtain a collection $\cC_1$ of $n_i-4$ cycles in~$G_i$.
Further, for each neighbor $z$ of $x$ in $G_{i-1}$, $(x,x'), (z,x'), (z,y'), (x',y')$ forms a cycle in $G_i$.
Applying this argument for all $N_{i-1}-(i-1)$ neighbors of $x$ in $G_{i-1}$, we obtain another collection $\cC_2$ of $N_{i-1}-(i-1)$ cycles in $G_i$.
Together, both collections form a set of $n_i - 4 + N_{i-1} - (i-1) = N_i - i-3$ cycles in $G_i$ such that $e$ is the only common edge of any two such cycles.

Secondly, we use a similar argument for the case $y = y'$.
By the inductive hypothesis, there are $N_{i-1}-(i-1)-3$ cycles of $G_{i-1}$ such that $\{x',y'\}$ is the only common edge of any two such cycles.
Further, $y$ has $n_i-1$ neighbors in~$\perm_{n_i}$.
Together, these induce $N_i-i-3$ cycles in $G_i$ that pairwise share $e$ as the only common edge.

The preceding two paragraphs complete the inductive proof.
This implies that $\perm$ is $(n-k-3)$-heavy.
Together with \cref{prop:lower_bound} and the fact that $\perm$ is $(n-k)$-regular, we then obtain the statement of the proposition.
\end{proof}

The main tool for the upper bound is the following proposition. 

\begin{proposition}
\label{prop:matchingscartesianproduct}
Let $G \eqdef H \bx B$ be the Cartesian product of a graph $H$ and a bipartite graph $B$. 
Suppose that $H$ has maximal matchings $N^\bullet$ and $N^\circ$ such that~$|N^\bullet | = |N^\circ|$ and every vertex of~$H$ is covered by $N^\bullet$ or $N^\circ$. 
Then~$G$ has maximal matchings $M^\bullet$ and $M^\circ$ such that $|M^\bullet| = |M^\circ| = |N^\bullet | \cdot |V(B)|$ and every vertex of $G$ is covered by $M^\bullet$ or $M^\circ$.
\end{proposition}

\begin{proof}
For $b \in V(B)$, let $H_b$ be the subgraph of $G$ induced by~${V(H) \times \{b\} \subseteq V(G)}$.
Let $N_b^\bullet$ and $N_{b}^\circ$ be maximal matchings in $H_b$ corresponding to~$N^\bullet$ and~$N^\circ$, respectively.
Let $B^\bullet$ and $B^\circ$ be the parts of~$B$. Set
\[
M^\bullet \eqdef \bigcup_{b \in B^\bullet} N_b^\bullet \cup \bigcup_{b \in B^\circ} N_b^\circ,
\qquad\text{and}\qquad
M^\circ \eqdef \bigcup_{b \in B^\bullet} N_b^\circ \cup \bigcup_{b \in B^\circ} N_b^\bullet.
\]
Clearly, $M^\bullet$ and $M^\circ$ are matchings in $G$, as they only use matching edges inside~$H_b$. 
For the cardinalities, we have
\[
|M^\bullet| = |N^\bullet| \cdot |B^\bullet| + |N^\circ| \cdot |B^\circ| = |N^\bullet | \cdot |V(B)| = |N^\circ| \cdot |B^\bullet| + |N^\bullet| \cdot |B^\circ| = |M^\circ|.
\]
Moreover, $M^\bullet$ is maximal since the exposed vertices of~$N_b^\bullet$ and~$N_{b'}^\circ$ are disjoint if~$bb' \in E(B)$ (and similarly, $M^\circ$ is maximal).
Finally, $M^\bullet \cup M^\circ$ covers~$V(G)$ since $N_b^\bullet \cup N_b^\circ$ covers~$V(H)$ for each~$b \in B$.
\end{proof}

\begin{remark}
Note that \cref{prop:matchingscartesianproduct} extends straightforward to a Cartesian product~$G \eqdef  H \bx K$ of two graphs~$H$ and~$K$ such that there are maximal matchings~$M_1, \dots, M_k$ of~$G$ and a coloring~$f : V(K) \to [k]$ such that~$M_{f(u)} \cup M_{f(v)}$ covers~$V(H)$ for any edge~$uv$ of~$K$.
\end{remark}

The upper bound in \cref{thm:quot_tight} then follows from \cref{ex:ub_hypercube} and the following proposition.

\begin{proposition}
\label{prop:ub_cartesian}
Let $n_1, \dots, n_k \geq 2$ be integers with $n_1 \geq 3$. Then the Cartesian product $\perm \eqdef \perm_{n_1} \bx \dots \bx \perm_{n_k}$ has a maximal matching of size $|V(\perm)| /3$. 
\end{proposition}

\begin{proof}
We show by induction on $i \in [k]$ that $G_i \eqdef \perm_{n_1} \bx \dots \bx \perm_{n_i}$ has maximal matchings $M_i^\bullet$ and $M_i^\circ$ of size $|V(G_i)|/3$ such that every vertex in $G$ is covered by $M_i^\bullet$ or $M_i^\circ$. 
For $i = 1$, this follows from Theorem~\ref{thm:ub_permutahedron} and Remark~\ref{rem:secondmatchingpn}.
For the induction step, assume that the above statement holds for $G_{i-1}$. 
By Proposition~\ref{prop:matchingscartesianproduct}, using that $\perm_{n_i}$ is bipartite, we obtain maximal matchings $M_i^\bullet$ and $M_i^\circ$ of size $|V(G_{i})|/3$ in $G_i$ such that every vertex in $G_i$ is covered by $M_i^\bullet$ or $M_i^\circ$. In particular, $\perm = G_k$ has a maximal matching of size $|V(\perm)|/3$, which proves the claim.
\end{proof}


\section{Open questions}

We conclude with a few open questions.
\begin{itemize}
	\item Is there a simple or constructive proof for conditions (ii) and (iii) of \cref{prop:ub_Hall} for the permutahedron?
	Is there a symmetric chain decomposition in the permutahedron?
	See \cref{rem:scd}.
	\item Can we get upper and lower bounds on the size of a minimum maximal matching for larger classes of polytopes generalizing the permutahedron and associahedron, in particular quotientopes (even only bipartite quotientopes or regular quotientopes), graphical zonotopes, graph associahedra, or Coxeter permutahedra?
	\item What can be said about minimal maximal matchings in Cartesian products in general? 
\end{itemize}

\bibliographystyle{plain}
\bibliography{refs}
\end{document}